\documentclass{amsart}
\title[On alternating sums of binomial...]{On alternating sums of binomial and $q$-binomial coefficients}
\usepackage{amssymb,amsmath,amsthm,epsfig,graphics,latexsym}
\usepackage{enumerate}
 \theoremstyle{definition}
 \newtheorem{definition}{Definition}
 \numberwithin{equation}{section}
  \theoremstyle{plain}

  \newtheorem{theorem}    {Theorem}

  \theoremstyle{remark}

  \DeclareMathOperator{\sgn}{sgn}
  \DeclareMathOperator{\inv}{inv}
\begin{document}
  \author{Mohamed El Bachraoui}
  %\address{Dept. Math. Sci,
 %United Arab Emirates University, PO Box 17551, Al-Ain, UAE}
 \email{melbachraoui@uaeu.ac.ae}
  \begin{abstract}
  In this paper we shall evaluate two alternating sums of binomial coefficients by a combinatorial argument. Moreover, by combining the same combinatorial idea with partition theoretic techniques, we provide $q$-analogues involving the $q$-binomial coefficients.
   \end{abstract}
  %\date{\textit{\today}}
  \maketitle
\section{Introduction}
Recall that the $q$-shifted factorials are given by
\[
(a;q)_0= 1,\quad (a;q)_n = \prod_{i=0}^{n-1}(1-a q^i),\quad
(a;q)_{\infty} = \lim_{n\to\infty}(a;q)_n =\prod_{i=0}^{\infty}(1-a q^i)
\]
and the $q$-binomial coefficients are given by
\[
{n\brack m} = \begin{cases}
\frac{(q;q)_n}{(q;q)_m (q;q)_{n-m}} , & \text{if\ } n\geq m\geq 0,\\
0, & \text{otherwise.}
\end{cases}
\]
Evaluating alternating sums and differences involving the binomial coefficients and finding their $q$-analogues involving the $q$-binomial coefficients have been extensively studied throughout the years and there is a rich literature
on the topic, see for instance \cite{Andrews-et-Al, Bailey, Bressoud-1, Bressoud-2, Gould, Guo-Zhang, Ismail-Kim-Stanton, Jackson}. A special case of a result by Andrews~ et~ al ~\cite{Andrews-et-Al} states that
\[
\sum_{k} (-1)^k {m+n\choose m-kl} \geq 0 \quad \text{if\ } |m-n|\leq l
\]
with a corresponding $q$-analogue stating that
\[
\sum_{k} (-1)^k q^{\frac{k^2 l(a+b)+kl(b-a)}{2}} {m+n \brack m-kl}
\]
is a polynomial in $q$ with nonnegative coefficients where $m,\ n,\ l,\ a,\ b$ are nonnegative integers
such that $a+b<2l$ and $b-l\leq n-m \leq l-a$. The authors proved their results using integer partitions.
Ismail~ et~ al \cite{Ismail-Kim-Stanton} extended the previous results by considering among other things expressions of the form
\[
\sum_{k} {m+n \choose m-kl} \cos(kx)
\]
and their $q$-analogues.
Recently Guo~and~Zhang~\cite{Guo-Zhang} gave combinatorial proofs for a variety of alternating
sums and differences of binomial and $q$-binomial coefficients including
\begin{equation}\label{mod2}
\sum_{k=-\infty}^{\infty}(-1)^k {2n \choose n+2k}=2^n
\end{equation}
and
\begin{equation}\label{mod3}
\sum_{k=-\infty}^{\infty}(-1)^k {2n \choose n+3k}=\begin{cases}
2 \cdot 3^{n-1}, & \text{if\ } n\geq 1,\\
1, & \text{if\ } n=0.
\end{cases}
\end{equation}
and their $q$-analogues
\begin{equation}\label{q-mod2-1}
\sum_{k=-\infty}^{\infty}(-1)^k q^{2k^2}{2n \brack n+2k}= (-q;q^2)_n
\end{equation}
and
\begin{equation}\label{q-mod3-1}
\sum_{k=-\infty}^{\infty}(-1)^k q^{(9k^2+3k)/2}{2n \brack n+3k}=
\begin{cases}
1, & \text{if\ } n=0, \\
(1+q^n)\frac{(q^3;q^3)_{n-1}}{(q;q)_{n-1}}, & \text{if\ } n\geq 1
\end{cases}
\end{equation}
respectively.
In this note we will prove generalisations of the identities (\ref{mod2})~and~(\ref{mod3}) along with generalisations of their
$q$-analogues~(\ref{q-mod2-1})~and~ (\ref{q-mod3-1}). Our $q$-analogues are expressed in terms of
certain restricted integer
partitions which we introduce now.
\begin{definition}
For any nonnegative integers $N$ and $M$  and any positive rational numbers $a\leq b$ let
$p_d(N, [a,b],M)$ denote the number of partitions of $N$ into exactly $M$ distinct parts all of which are in the integer interval $[a,b]$. We assume that $p_d(-N,[a,b],M)=0$ if $N\leq 0$.
\end{definition}
\noindent
Further, to simplify the formulas, we introduce the following notation.
\begin{definition}
Let $m$ and $n$ be nonnegative integers and let
$a_{j,l}=p_d(j,[n,n+m],l)$. For any positive integer $k$ and any nonnegative integer $r$ we
 let $A_{n,m,k,r}(q)$ be the polynomial in $q$ given by:
\[
A_{n,m,k,r}(q) = A_{k,r}(q)
\]
\[=
\sum_{l=0}^{\lfloor\frac{m+1-r}{k}\rfloor} \sum_{j=(kl+r)(n+\frac{kl-1+r}{2})}^{(kl+r)(n+m-\frac{kl-1+r}{2})}
a_{j,kl+r}q^j +
\sum_{l=0}^{\lfloor\frac{m+2-r}{k}\rfloor}
 \sum_{j=(kl+r-1)(n+\frac{kl-2+r}{2})}^{(kl+r-1)(n+m-\frac{kl-2+r}{2})}
a_{j,kl+r-1}q^j
\]
\end{definition}
\section{The results}
\begin{theorem}\label{MainSum}
If $m$ is an integer and $n$ is a nonnegative integer, then
\[
\sum_{k=-\infty}^{\infty} (-1)^k {2n+m\choose n+2k} = 2^{n+\frac{m}{2}}\cos\frac{m\pi}{4}
\]
\end{theorem}
\begin{theorem}\label{q-analog-1}
If $m$ is a nonnegative integer and $n$ is a positive integer, then
\begin{align*}
\sum_{k=-\infty}^{\infty} (-1)^k  q^{2 k^2}{2n+m\brack n+2k} &=
(-q;q^2)_n  \\
& \times \left(\sum_{l=0}^{\lfloor\frac{m}{4}\rfloor}
 \sum_{j=2l(2n+4l)}^{2l(2n+2m-4l)} a_{j,l} q^{j}
   -
\sum_{l=0}^{\lfloor\frac{m-2}{4}\rfloor}
 \sum_{j=(2l+1)(2n+4l+2)}^{(2l+1)(2n+2m-4l-2)} b_{j,l} q^{j}  \right),
\end{align*}
where
\[
a_{j,l}= p_d(j,[n+1,n+m-1],4l)\ \text{and\ }
b_{j,l}=p_d(j,[n+1,n+m-1],4l+2).
\]
\end{theorem}
\begin{theorem}\label{q-analog-2}
If $n$ and $m$ are positive integers, then
\begin{align*}
\sum_{k=-\infty}^{\infty} (-1)^k  q^{2 k^2-k}{2n+m\brack n+2k}  &=
(-q^2;q^2)_n  \\
& \times \left\{ \sum_{l=0}^{\lfloor\frac{m-1}{4}\rfloor}
 \sum_{j=2l(2n+4l+1)}^{2l(2n+2m-4l-1)} a_{j,l}q^j   + \sum_{l=0}^{\lfloor\frac{m}{4}\rfloor}
   \sum_{j=(4l-1)(n+2l)}^{(4l-1)(n+m-2l)} c_{j,l} q^j   \right. \\
& \left.\vphantom{\int_t}
  -
\left(  \sum_{l=0}^{\lfloor\frac{m-3}{4}\rfloor}
  \sum_{j=(2l+1)(2n+4l+3)}^{(2l+1)(2n+2m-4l-3)} b_{j,l} q^j +
  \sum_{l=0}^{\lfloor\frac{m-2}{4}\rfloor}\sum_{j=(4l+1)(n+2l+1)}^{(4l+1)(n+m-2l-1)} d_{j,l}q^j    \right)
  \right\},
\end{align*}
where
\[
a_{j,l}=p_d(j,[n+1,n+m-1],4l),\
c_{j,l}=p_d(j,[n+1,n+m-1],4l-1),
\]
\[
b_{j,l}=p_d(j,[n+1,n+m-1],4l+2),\
d_{j,l}=p_d(j,[n+1,n+m-1],4l+1).
\]
\end{theorem}
 \begin{theorem}\label{MainSum-2}
If $m$ is an integer and $n$ is a nonnegative integer, then
\[
\sum_{k=-\infty}^{\infty} (-1)^k {2n+m\choose n+3k} = \begin{cases}
1, & \text{if\ } n=m=0, \\
2\cdot 3^{n-1+\frac{m}{2}}\cos\frac{m\pi}{6}, & \text{otherwise.}
\end{cases}
\]
\end{theorem}
\begin{theorem}\label{q-analog-2-1}
Let $m$ be a nonnegative integer and let $n$ be a positive integer. Then
\[
\sum_{k=-\infty}^{\infty} (-1)^k  q^{\frac{9 k^2 + 3k}{2}}{2n+m\brack n+3k}  =
\frac{(q^3;q^3)_{n-1}}{(q;q)_{n-1}} \Bigl(A_{6,0}(q)-A_{6,3}(q)+A_{6,1}(q)-A_{6,4}(q) \Bigr).
\]
\end{theorem}
\noindent
{\bf Remarks.}
 {\bf 1.} As a consequence of the previous theorems, sums involving the partitions $p_d(N,[a,b],M)$
 and formulas involving $A_{n,m,6,r}(1)$
will be obtained upon letting $q\to 1$.
For instance, letting $m=4N+2$ and $q\to 1$ we find by combining Theorem~\ref{MainSum} with
Theorem~\ref{q-analog-1} that
\[
\sum_{l=0}^{\lfloor\frac{m}{4}\rfloor}
 \sum_{j=2l(2n+4l)}^{2l(2n+2m-4l)} p_d(j,[n+1,n+m-1],4l) =
 \]
 \[
\sum_{l=0}^{\lfloor\frac{m-2}{4}\rfloor}
 \sum_{j=(2l+1)(2n+4l+2)}^{(2l+1)(2n+2m-4l-2)} p_d(j,[n+1,n+m-1],4l+2)
 \]
 and letting $m=6N+3$ and $q\to 1$ we find by combining Theorem~\ref{MainSum-2} with
Theorem~\ref{q-analog-2-1} that for all positive integer $n$
\[
A_{n,m,6,0}(1)-A_{n,m,6,3}(1)+A_{n,m,6,1}(1)-A_{n,m,6,4}(1) = 0.
 \]
{\bf 2.} Notice that for any integers $m$ and $n$, the binomial and $q$-binomial coefficients make the series
 in the previous theorems finite on both ends. However, if $m, n \to +\infty$, then
 by virtue of the Jacobi triple product (see~\cite{Andrews-Askey-Roy, Gasper-Rahma}) the sum in Theorem~\ref{q-analog-1} becomes:
\[
\lim_{m,n\to\infty}\sum_{k=-\infty}^{\infty} (-1)^k  q^{2 k^2}{2n+m\brack n+2k} =
\frac{(q^2;q^4)_{\infty}^2 (q^4;q^4)_{\infty}}{(q;q)_{\infty}} = (-q;q)_{\infty}(q^2;q^4)_{\infty}
\]
and similarly the sums in Theorem~\ref{q-analog-2} and
Theorem~\ref{q-analog-2-1} respectively become:
\[
\lim_{m,n\to\infty}\sum_{k=-\infty}^{\infty} (-1)^k  q^{2 k^2-k}{2n+m\brack n+2k} =
\frac{1}{(q^2;q^4)_{\infty}}
\]
and
\[
\lim_{m,n\to\infty}\sum_{k=-\infty}^{\infty} (-1)^k  q^{\frac{9 k^2 +3 k}{2}}{2n+m\brack n+3k} =
\frac{(q^3;q^3)_{\infty}}{(q;q)_{\infty}}.
\]
{\bf 3.} We note further that the sums in Theorems \ref{q-analog-1}, \ref{q-analog-2}, and \ref{q-analog-2-1} are related  to a finite sum version of a $_{2}\psi_{2}$ sum.
Refer to \cite{Andrews-Askey-Roy, Gasper-Rahma} for more details about the function $_{2}\psi_{2}$.
%\end{rmks}
%
%Throughout let $\N$  be the set of positive integers and let $\N_0$ be the set of nonnegative integers.
%
\section{Proof of Theorem~\ref{MainSum}}
Suppose first that $n>0$ and $m\geq 0$.
Following Guo~and~Zhang~\cite{Guo-Zhang},
throughout  let
$S= \{ a_1,\ldots,a_{2n}, a_{2n+1},\ldots, a_{2n+m}\}$ be a set of $2n+m$ elements.
Let
\[
\begin{split}
\mathcal{F}&= \{ A\subseteq S:\ \# A\equiv n \pmod{2} \},
\\
\mathcal{G} &= \{ A\subseteq S:\ \#(A\cap\{a_{2i-1},a_{2i}\})=1\ \text{for all\ } i=1,\ldots, n\},
\\
\mathcal{H}&= \{A\in F:\ \#(A\cap\{a_{2i-1},a_{2i}\})\not=1\ \text{for some\ } i=1,\ldots, n\}.
\end{split}
\]
For simplicity of notation, if $A\in\mathcal{S}$, we let
\[
A' = A\cap \{a_{2n+1},\ldots,a_{2n+m} \}.
\]
We define a map
\[
\sgn: \mathcal{F}\to \{-1,1\},\quad
\sgn(A) = (-1)^{\frac{\# A-n}{2}}.
\]
Then it is clear that
\begin{equation}\label{SumonSgn}
\sum_{k=-\lfloor\frac{n}{2}\rfloor}^{\lfloor\frac{n+m}{2}\rfloor}(-1)^k {2n+m\choose n+2k} =
\sum_{A\in\mathcal{F}} \sgn(A).
\end{equation}
Furthermore, if $A\in  \mathcal{H}$ let $i_A$ be the minimum index $i$ such that
$\#(A\cap\{a_{2i-1},a_{2i} \}) \not= 1$. Next we define a map
\[
\inv: \mathcal{H} \to \mathcal{H}
\]
\[
\inv(A) = \begin{cases}
A\cup\{a_{2i_A -1}, a_{2i_A} \} &  \text{if\ } \#(A\cap\{a_{2i_A-1},a_{2i_A}\}) = 0\\
A\setminus\{a_{2i_A -1}, a_{2i_A} \} &  \text{if\ } \#(A\cap\{a_{2i_A -1},a_{2i_A}\}) = 2.
\end{cases}
\]
Then it is easily seen that the function $\text{inv}$ is an involution satisfying
$\sgn(\inv(A)) = - \sgn(A)$ and therefore we have
\begin{align}\label{SplitSgn}
\sum_{A\in\mathcal{F}} \sgn(A) &=
\sum_{\substack{A\in\mathcal{G}\\ \#A' \equiv 0\pmod{2}}} \sgn(A)
\nonumber \\
&=\sum_{\substack{A\in\mathcal{G}\\ \#A' \equiv 0\pmod{4}}} 1 -
\sum_{\substack{A\in\mathcal{G}\\ \#A' \equiv 2\pmod{4}}} 1
\nonumber \\
&= 2^{n}\left(\sum_{l=0}^{\lfloor\frac{m}{4}\rfloor} {m\choose 4l}-
\sum_{l=0}^{\lfloor\frac{m-2}{4}\rfloor} {m\choose 2+4l} \right) .
\end{align}
Then the case $m=0$ follows easily  from  the identities (\ref{SumonSgn}) and (\ref{SplitSgn}).
 Suppose now that $m\geq 1$. Note that by the formulas (\ref{SumonSgn}) and (\ref{SplitSgn}) it suffices to show that
\[
\sum_{l\geq 0} {m\choose 4l}  - \sum_{l\geq 0} {m\choose 4l+2} = 2^{\frac{m}{2}} \cos\frac{m\pi}{4}.
\]
By a well-known result, see Gould~\cite{Gould} and Benjamin et al \cite{Benjamin-et-Al}, we have
\[
\sum_{l\geq 0} {m\choose 4l} =
\frac{1}{4}\sum_{k=0}^3 (1+ e^{i\frac{2\pi k}{4}})^m
\]
and
\[
\sum_{l\geq 0} {m\choose 4l+2}  = \frac{1}{4}\sum_{k=0}^3 e^{-2k \frac{i2\pi}{4}}
(1+ e^{i\frac{2\pi k}{4}})^m,
\]
from which it follows that
\[
\sum_{l\geq 0} {m\choose 4l}  - \sum_{l\geq 0} {m\choose 4l+2} = 2^{\frac{m}{2}} \cos\frac{m\pi}{4},
\]
as desired. Suppose next that $m<0$ and let $M=-m = 4s-r$ with $0\leq r <4$. Then with the help of the previous case we have
\[
\begin{split}
\sum_{k=-\infty}^{\infty} (-1)^k{2n+m\choose n+2k} &=
\sum_{k=-\infty}^{\infty} (-1)^k {2n-4s +r\choose n-2s+2(k+s)} \\
&=
(-1)^s\sum_{k=-\infty}^{\infty}(-1)^k {2(n-2s)+r \choose (n-2s)+2k} \\
&=
(-1)^s 2^{n-2s+\frac{r}{2}} \cos\frac{\pi r}{4} \\
&= 2^{n-\frac{M}{2}} \cos\frac{\pi M}{4},
\end{split}
\]
as desired. Further, suppose that $n=0$ and $m\geq 4$. Then
\[
\begin{split}
\sum_{k=-\infty}^{\infty} (-1)^k {m\choose 2k} &=
\sum_{k=-\infty}^{\infty} (-1)^k {4+(m-4)\choose 2+2(k-1)} \\
 &=
-\sum_{k=-\infty}^{\infty} (-1)^k {4+(m-4)\choose 2+2k} \\
&=
-2^{2+\frac{m-4}{2}} \cos\frac{(m-4)\pi}{4} \\
&=
2^{\frac{m}{2}}\cos\frac{m\pi}{4}.
\end{split}
\]
Finally it is easy to check the cases $n=0$ and $m=1, 2, 3$. This completes the proof.
\section{Proof of Theorem~\ref{q-analog-1}}
If the $2n+m$ elements of the set $S$ are integers and $A\subseteq S$ then we define the \emph{weight} of $A$ by $\|A\|=\sum_{a\in A} a$. From the well-known fact, see Andrews~\cite{Andrews-1},
\[
(zq;q)_N = \sum_{j=0}^N {N\brack j}(-1)^j z^j q^{{j+1\choose 2}}
\]
we conclude that
\begin{equation}\label{weight-sum}
\sum_{\substack{A\subseteq\{1,\ldots,n\}\\ \#A=k}}q^{\|A\|} = {n\brack k} q^{{k+1\choose 2}}.
\end{equation}
Now let
\[
S=\{i-\frac{2n+1}{2}:\ i=1,\ldots,2n+m\} = \{\pm\frac{1}{2},\pm\frac{3}{2},\ldots,\pm\frac{2n-1}{2},
\frac{2n+1}{2},\ldots, \frac{2n+2m-1}{2}\},
\]
that is,
\[ a_{2i-1}= - \frac{2i-1}{2}, \ a_{2i} = \frac{2i-1}{2} \ \text{for\ }
i=1,\ldots,n \ \text{and\ }
 a_{2n+j} = \frac{2n+2j-1}{2}\ \text{for\ } j=1,\ldots,m.
 \]
Then the function $\inv$ defined above is weight preserving and therefore we have
\begin{equation}\label{weight-preserving}
\sum_{A\in\mathcal{F}}\sgn(A) q^{\|A\|} =
 \sum_{\substack{A\in\mathcal{G}\\ \#A' \equiv 0\pmod{4}}} q^{\|A\|} -
\sum_{\substack{A\in\mathcal{G}\\ \#A' \equiv 2\pmod{4}}} q^{\|A\|}
\end{equation}
As to the left-hand-side of the relation~(\ref{weight-preserving}), we have
\[
\begin{split}
\sum_{A\in\mathcal{F}}\sgn(A) q^{\|A\|}
&=
\sum_{k=-\lfloor\frac{n}{2}\rfloor}^{\lfloor\frac{n+m}{2}\rfloor}
 \sum_{\substack{A\subseteq S\\ \#A=n+2k}} \sgn(A) q^{\|A\|}  \\
&=
 \sum_{k=-\lfloor\frac{n}{2}\rfloor}^{\lfloor\frac{n+m}{2}\rfloor}
 \sum_{\substack{A\subseteq \{1,\ldots,2n+m\} \\ \#A=n+2k}} (-1)^k q^{\|A\|+(n+2k)(-\frac{2n+1}{2})} \\
 &=
 \sum_{k=-\lfloor\frac{n}{2}\rfloor}^{\lfloor\frac{n+m}{2}\rfloor}
 (-1)^k q^{-\frac{2n+1}{2}(n+2k)}  \sum_{\substack{A\subseteq \{1,\ldots,2n+m\} \\ \#A=n+2k}} q^{\|A\|} ,
 \end{split}
 \]
 which with the help of the identity (\ref{weight-sum}) gives
 \begin{equation}\label{lhs}
 \begin{split}
 \sum_{A\in\mathcal{F}}\sgn(A) q^{\|A\|}
 &=
 \sum_{k=-\lfloor\frac{n}{2}\rfloor}^{\lfloor\frac{n+m}{2}\rfloor}
 (-1)^k q^{-\frac{2n+1}{2}(n+2k)} {2n+m\brack n+2k} q^{ {n+2k+1\choose 2}} \\
 &=
 \sum_{k=-\lfloor\frac{n}{2}\rfloor}^{\lfloor\frac{n+m}{2}\rfloor}
 (-1)^k q^{\frac{(2k)^2 - n^2}{2} }{2n+m\brack n+2k}.
  \end{split}
 \end{equation}
 As to the first sum on the right-hand-side of the relation~(\ref{weight-preserving}), notice first
 that the least sum and the largest sum into exactly $4l$ distinct parts belonging to the set
$\{(2n+1)/2,\ldots,(2n+2m-1)/2\}$ are respectively
\[
\frac{2n+1}{2}+\ldots+ \frac{2n+2(4l-1)+1}{2} =2l(2n+4l)
\]
and
\[
\frac{2n+2(m-4l)+1}{2}+\ldots+\frac{2n+2m-1}{2} = 2l(2n+2m-4l).
\]
 So, letting
 \[
 a_{j,l}=p_d \left( j,\left[\frac{2n+1}{2},\frac{2n+2m-1}{2}\right],4l \right) = p_d(j,[n+1,n+m-1],4l ), \\
 \]
 it is easily checked that
 \begin{align}\label{rhs1}
 \sum_{\substack{A\in\mathcal{G}\\ \#A' \equiv 0\pmod{4}}} 1
&=
\prod_{i=1}^n (q^{-\frac{2i-1}{2}}+ q^{\frac{2i-1}{2}}) \times
\sum_{l=0}^{\lfloor\frac{m}{4}\rfloor}\sum_{j=2l(2n+4l)}^{2l(2n+2m-4l)}
a_{j,l} q^{j}  \nonumber \\
&=
q^{-\frac{n^2}{2}}  (-q;q^2)_n \sum_{l=0}^{\lfloor\frac{m}{4}\rfloor}\sum_{j=2l(2n+4l)}^{2l(2n+2m-4l)}
a_{j,l} q^{j}.
\end{align}
As to the second sum on the right-hand-side of the relation~(\ref{weight-preserving}) we
use the same remark as before with $4l$ replaced by $4l+2$ and $a_{j,l}$
replaced by
\[
b_{j,l}=p_d \left( j,\left[\frac{2n+1}{2},\frac{2n+2m-1}{2}\right],4l+2 \right) = p_d(j,[n+1,n+m-1],4l +2)
\]
to obtain
\begin{align}\label{rhs2}
\sum_{\substack{A\in\mathcal{G}\\ \#A' \equiv 2\pmod{4}}}  1
&=
\prod_{i=1}^n (q^{-\frac{2i-1}{2}}+ q^{\frac{2i-1}{2}}) \times
\sum_{l=0}^{\lfloor\frac{m-2}{4}\rfloor}\sum_{j=(2l+1)(2n+4l+2)}^{(2l+1)(2n+2m-4l-2)}
b_{j,l} q^{j} \nonumber  \\
&=
q^{-\frac{n^2}{2}}  (-q;q^2)_n \sum_{l=0}^{\lfloor\frac{m-2}{4}\rfloor}
\sum_{j=(2l+1)(2n+4l+2)}^{(2l+1)(2n+2m-4l-2)}
b_{j,l} q^{j}.
\end{align}
Then the desired formula follows by combining the relation (\ref{weight-preserving}) with the relations (\ref{lhs}), (\ref{rhs1}), and  (\ref{rhs2}).
\section{Proof of Theorem~\ref{q-analog-2}}
\noindent
We proceed as in the proof of Theorem~\ref{q-analog-1}. Let
\[
S = \{i - (n+1):\ i=1,\ldots,2n+m\} = \{\pm 1,\ldots,\pm n,0,n+1,n+2,\ldots,n+m-1\},
\]
that is,
\[
a_{2i-1}=-i,\ a_{2i}=i\ \text{for \ } i=1,\ldots,n,\
a_{2n+1}=0,\ \text{and\ }
a_{2n+j} = n+j-1\ \text{for \ } j=2,\ldots,m.
\]
Then
\begin{equation}\label{lhs-1}
\begin{split}
\sum_{A\in\mathcal{F}}\sgn(A) q^{\|A\|}
& =
 \sum_{k=-\lfloor\frac{n}{2}\rfloor}^{\lfloor\frac{n+m}{2}\rfloor}
 (-1)^k q^{ (n+2k)(-n-1) } q^{ {n+2k+1 \choose 2} }{2n+m\brack n+2k} \\
& =
\sum_{k=-\lfloor\frac{n}{2}\rfloor}^{\lfloor\frac{n+m}{2}\rfloor}
 (-1)^k q^{2k^2- k - \frac{n(n+1)}{2}} {n+m\brack n+2k}.
 \end{split}
\end{equation}
Let as before
\[
a_{j,l}=p_d(j,[n+1,n+m-1],4l),\  \text{and\ }
b_{j,l}=p_d(j,[n+1,n+m-1],4l+2).\
\]
However, because of the presence of $0$ among the elements of $\mathcal{S}$ we shall also take into
account partition into exactly $4l-1$ (nonzero) parts and therefore we let
\[
c_{j,l}=p_d(j,[n+1,n+m-1],4l-1)\ \text{and\ } d_{j,l}=p_d(j,[n+1,n+m-1],4l+1).
\]
We have
\begin{equation}\label{rhs1-1}
 \begin{split}
 \sum_{\substack{A\in\mathcal{G}\\ \#A' \equiv 0\pmod{4}}} q^{\|A\|}
&=
\prod_{i=1}^n (q^{-i}+ q^{i})   \\
& \times \sum_{l=0}^{\lfloor\frac{m-1}{4}\rfloor} \sum_{j=2l(2n+4l+1)}^{2l(2n+2m-4l-1)}
a_{j,l} q^j+
\sum_{l=0}^{\lfloor\frac{m}{4}\rfloor}\sum_{j=(4l-1)(n+2l)}^{(4l-1)(n+m-2l)} c_{j,l} q^j   \\
&=
q^{-\frac{n(n+1)}{2}}  (-q^2;q^2)_n \\
& \times \sum_{l=0}^{\lfloor\frac{m-1}{4}\rfloor}\left( \sum_{j=2l(2n+4l+1)}^{2l(2n+2m-4l-1)}
a_{j,l} q^j +  \sum_{j=(4l-1)(n+2l)}^{(4l-1)(n+m-2l)} c_{j,l} q^j \right) ,
\end{split}
\end{equation}
and
\[
\begin{split}
\sum_{\substack{A\in\mathcal{G}\\ \#A' \equiv 2\pmod{4}}} q^{\|A\|}
&=\prod_{i=1}^n (q^{-i}+ q^{i}) \\
& \times \sum_{l=0}^{\lfloor\frac{m-3}{4}\rfloor} \sum_{j=(2l+1)(2n+4l+3)}^{(2l+1)(2n+2m-4l-3)}
b_{j,l}  q^j +
\sum_{l=0}^{\lfloor\frac{m-2}{4}\rfloor}\sum_{j=(4l+1)(n+2l+1)}^{(4l+1)(n+m-2l-1)} d_{j,l} q^j
\end{split}
\]
\begin{equation}\label{rhs2-1}
\begin{split}
& \qquad\quad =
q^{-\frac{n(n+1)}{2}}  (-q^2;q^2)_n \\
& \qquad\quad \times \sum_{l=0}^{\lfloor\frac{m-3}{4}\rfloor}\left( \sum_{j=(2l+1)(2n+4l+3)}^{2l(2n+2m-4l-3)}
b_{j,l}  q^j +  \sum_{j=(4l+1)(n+2l+1)}^{(4l+1)(n+m-2l-1)} d_{j,l}  q^j \right).
\end{split}
\end{equation}
Then the desired formula follows by combining the relations (\ref{lhs-1}), (\ref{rhs1-1}), and
(\ref{rhs2-1}) with the relation (\ref{weight-preserving}).
\section{Proof of Theorem~\ref{MainSum-2}}
The case $n=m=0$ is clear. Suppose now that $n>0$ and $m\geq 0$.
Extending definitions from Guo~and~Zhang~\cite{Guo-Zhang},
throughout  let
$S= \{ a_1,\ldots,a_{2n}, a_{2n+1},\ldots, a_{2n+m}\}$ be a set of $2n+m$ elements and
let
\[
\begin{split}
\mathcal{F}&= \{ A\subseteq S:\ \# A\equiv n \pmod{3} \},
\\
\mathcal{G} &=
\{ A\subseteq \mathcal{F}:\ \#(A\cap\{a_{1},a_{2i+1}\})\not\in\{i-1,i+2\}\ \text{for all\ } i=1,\ldots, n-1\}.
%\\
%\mathcal{H}&= \{A\in F:\ \#(A\cap\{a_{2i-1},a_{2i}\})\not=1\ \text{for some\ } i=1,\ldots, n\}.
\end{split}
\]
We define a map
\[
\sgn: \mathcal{F}\to \{-1,1\},\quad
\sgn(A) = (-1)^{\frac{\# A-n}{3}}.
\]
Then it is clear that
\begin{equation}\label{SumonSgn}
\sum_{k=-\lfloor\frac{n}{3}\rfloor}^{\lfloor\frac{n+m}{3}\rfloor}(-1)^k {2n+m\choose n+3k} =
\sum_{A\in\mathcal{F}} \sgn(A).
\end{equation}
Furthermore, if $A\in  \mathcal{F}\setminus\mathcal{G}$ let $i_A$ be the minimum index $i$ such that
$\#(A\cap\{a_{1},a_{2i+1} \}) \in\{i-1,i+2\}$. Letting
\[
A' = A\Delta\{a_1,\ldots,a_{2i_A +1} \} = (A\cup\{a_1,\ldots,a_{2i_A +1} \} )\setminus
(A\cap\{a_1,\ldots,a_{2i_A +1} \} ),
\]
note that it is easily seen that $\# A'= \# A\pm 3$.
Next we define a map
\[
\inv: \mathcal{F}\setminus\mathcal{G} \to \mathcal{F}\setminus\mathcal{G},\quad
\inv(A) = A''
\]
as follows:
\begin{itemize}
\item $a_1\in A''$ if and only if $a\not\in A$,
\item $a_{2j}, a_{2j+1}\in A''$ if $a_{2j}, a_{2j+1}\not\in A$ ($j=1,\ldots, i_A$);
\item $a_{2j}, a_{2j+1}\not\in A''$ if $a_{2j}, a_{2j+1}\in A$ ($j=1,\ldots, i_A$);
\item $a_{2j}\in A''$ and $a_{2j+1}\not\in A''$, if $a_{2j}\in A$ and $a_{2j+1}\not\in A$ ($j=1,\ldots, i_A$);
\item $a_{2j}\not\in A''$ and $a_{2j+1}\in A''$, if $a_{2j}\not\in A$ and $a_{2j+1}\in A$ ($j=1,\ldots, i_A$);
\item $a_k \in A''$ if and only if $a_k \in A$ ($2i_A + 2 \leq k \leq 2n+m$).
\end{itemize}
Observing  that $\# A''= \# A'$, we have from the previous note that $\# A''=\# A \pm 3$. Further, it is clear
that $A''\in \mathcal{F}\setminus\mathcal{G}$ and that
 the function $\text{inv}$ is an involution satisfying
$\sgn(\inv(A)) = - \sgn(A)$. Therefore, we have
\begin{equation}\label{SplitSgn-1}
\sum_{A\in\mathcal{F}} \sgn(A) = \sum_{A\in\mathcal{G}} \sgn(A) +
\sum_{A\in\mathcal{F}\setminus\mathcal{G}} \sgn(A) = \sum_{A\in\mathcal{G}} \sgn(A).
\end{equation}

\noindent
We now evaluate $\sum_{A\in\mathcal{G}} \sgn(A)$. We claim that if $A\in\mathcal{G}$, then
\[
\# (A\cap\{a_1,\ldots,a_{2i+1}\}) \in\{i,i+1\},\ \text{for all\ } i=1,\ldots,n-1.
\]
Indeed, for $i=1$ the claim is evident. Next suppose that the statement holds for $i-1$. Then
\[
\# (A\cap\{a_1,\ldots,a_{2i+1}\}) \in\{ i-1,i,i+1,i+2\},
\]
from which the claim follows since the cases $i-1$ and $i+2$ are excluded by definition.
For simplicity of notation, if $A\in\mathcal{S}$, we let
\[
A_1= A\cap \{a_{1},\ldots,a_{2n-1} \},\quad
A_2 = A\cap \{a_{2n},\ldots,a_{2n+m} \},
\]
\[
\mathcal{G}_1 = \{A\in\mathcal{G}: \# A_1 = n\},\ \text{and\ }
\mathcal{G}_2 = \{A\in\mathcal{G}: \# A_1 = n-1\}.
\]
By the previous claim we have $\# A_1\in\{ n-1,n\}$,
which combined with the identity~(\ref{SplitSgn-1}) yields
\begin{equation}\label{SplitSgn-2}
\sum_{A\in\mathcal{F}} \sgn(A) = \sum_{A\in\mathcal{G}_1} \sgn(A) +
\sum_{A\in\mathcal{G}_2} \sgn(A).
\end{equation}
Furthermore, we clearly have
\begin{equation}\label{SubSplit-1}
\sum_{A\in\mathcal{G}_1} \sgn(A) =
\sum_{\substack{A\in\mathcal{G}_1 \\ \# A_2 \equiv 0 \pmod{6}}}1 -
\sum_{\substack{A\in\mathcal{G}_1 \\ \# A_2 \equiv 3 \pmod{6}}}1
\end{equation}
and
\begin{equation}\label{SubSplit-2}
\sum_{A\in\mathcal{G}_2} \sgn(A) =
\sum_{\substack{A\in\mathcal{G}_2\\   \# A_2 \equiv 1 \pmod{6}}}1 -
\sum_{\substack{A\in\mathcal{G}_2\\  \# A_2 \equiv 4 \pmod{6}}}1.
\end{equation}
Moreover, if $A\in\mathcal{G}$, then there are three possible choices for
$A \cap\{a_{2i},a_{2i+1}\}$ for all $i=1,\ldots,n-1$ and there is exactly one possible choice for
$A \cap\{a_1\}$, implying that
\begin{equation}\label{G-cardinality}
\# \mathcal{G}_1= \# \mathcal{G}_2 = 3^{n-1}.
\end{equation}
Then from the relations~(\ref{SubSplit-1}) and
(\ref{G-cardinality}) we derive
\begin{equation}\label{SubSplit-3}
\sum_{A\in\mathcal{G}_1} \sgn(A) = 3^{n-1}\left(
\sum_{l=0}^{\lfloor\frac{m+1}{6}\rfloor}{m+1\choose 6l} -
\sum_{l=0}^{\lfloor\frac{m-2}{6}\rfloor}{m+1\choose 6l+3} \right)
\end{equation}
and from the relations~(\ref{SubSplit-2}) and (\ref{G-cardinality}) we derive
\begin{equation}\label{SubSplit-4}
\sum_{A\in\mathcal{G}_2} \sgn(A) = 3^{n-1}\left(
\sum_{l=0}^{\lfloor\frac{m}{6}\rfloor}{m+1\choose 6l+1} -
\sum_{l=0}^{\lfloor\frac{m-3}{6}\rfloor}{m+1\choose 6l+4} \right).
\end{equation}
Substituting (\ref{SubSplit-3}) and (\ref{SubSplit-4}) in the formula (\ref{SplitSgn-2}) gives
\begin{equation}\label{SplitSgn-3}
\begin{split}
\sum_{A\in\mathcal{F}} \sgn(A)
&= 3^{n-1}\left(
\sum_{l=0}^{\lfloor\frac{m+1}{6}\rfloor}{m+1\choose 6l} -
\sum_{l=0}^{\lfloor\frac{m-2}{6}\rfloor}{m+1\choose 6l+3}  \right) \\
& \quad +3^{n-1}\left(
\sum_{l=0}^{\lfloor\frac{m}{6}\rfloor}{m+1\choose 6l+1} -
\sum_{l=0}^{\lfloor\frac{m-3}{6}\rfloor}{m+1\choose 6l+4} \right).
\end{split}
\end{equation}
Then the case $m=0$ follows easily  from  the identities (\ref{SumonSgn}) and (\ref{SplitSgn-3}).
 Suppose now that $m\geq 1$.
 By a well-known result, see Gould~\cite{Gould} and Benjamin et al \cite{Benjamin-et-Al}, we have
\begin{equation}\label{GouldIdentity}
\sum_{l\geq 0}{N\choose a+rl} = \frac{1}{r}\sum_{j=0}^{r-1}
e^{-i\frac{2\pi j a}{r}} (1+e^{-i\frac{2\pi j a}{r}})^N,
\end{equation}
where $r, a$, and $N$ are nonnegative integers such that $a<r$. An application of this relation along with basic simplifications imply
\[
\sum_{l\geq 0} {m\choose 6l}  - \sum_{l\geq 0} {m\choose 6l+3} =
2\cdot 3^{\frac{m-1}{2}} \cos\frac{(m+1)\pi}{6}
\]
and
\[
\sum_{l\geq 0} {m\choose 6l+1}  - \sum_{l\geq 0} {m\choose 6l+4} =
2\cdot 3^{\frac{m-1}{2}} \cos\frac{(m-1)\pi}{6},
\]
which combined with ((\ref{SplitSgn-3})) yields
\begin{equation}\label{SplitSgn-4}
\sum_{A\in\mathcal{F}} \sgn(A) = 3^{n-1} (2\cdot 3^{\frac{m}{2}} \cos\frac{m\pi}{6}) =
2\cdot 3^{n-1+\frac{m}{2}} \cos\frac{m\pi}{6},
\end{equation}
as desired. Suppose next that $m<0$ and let $M=-m = 6s-r$ with $0\leq r <6$. Then with the help of the above we have
\[
\begin{split}
\sum_{k=-\infty}^{\infty} {2n+m\choose n+3k} &=
\sum_{k=-\infty}^{\infty}{2n-6s +r\choose n-3s+3(k+s)} \\
&=
(-1)^s\sum_{k=-\infty}^{\infty}(-1)^k {2(n-3s)+r \choose (n-3s)+3k} \\
&=
(-1)^s 2\cdot 3^{n-3s-1+\frac{r}{2}} \cos\frac{\pi r}{6} \\
&= 2\cdot 3^{n-1-\frac{M}{2}} \cos\frac{\pi M}{6},
\end{split}
\]
as desired. Further, suppose that $n=0$ and $m\geq 6$. Then
\[
\begin{split}
\sum_{k=-\infty}^{\infty} {m\choose 3k} &=
\sum_{k=-\infty}^{\infty} {6+(m-6)\choose 3+3(k-1)} \\
 &=
-\sum_{k=-\infty}^{\infty} {6+(m-6)\choose 3+3k} \\
&=
-2\cdot 3^{3-1+\frac{m-6}{2}} \cos\frac{(m-6)\pi}{6} \\
&=
2\cdot 3^{-1+\frac{m}{2}}\cos\frac{m\pi}{6}.
\end{split}
\]
Finally it is easy to check the cases $n=0$ and $m\in\{1, 2, 3,4,5\}$. This completes the proof.
\section{Proof of Theorem ~\ref{q-analog-2-1}}
Let
\[
S=\{i-n:\ i=1,\ldots,2n+m\} = \{\pm 1,\ldots,\pm(n-1), 0, n, n+1,\ldots,n+m\},
\]
with
\[ a_1 = 0, a_{2i}= - i, \ a_{2i+1} = i \ \text{for\ }
i=1,\ldots,n-1 \ \text{and\ }
 a_{2n+j} = n+j \ \text{for\ } j=0,\ldots,m.
 \]
Then the function $\inv$ defined above is weight preserving and therefore we have
\begin{equation}\label{weight-preserving}
\begin{split}
\sum_{A\in\mathcal{F}}\sgn(A) q^{\|A\|} &=
\sum_{A\in\mathcal{G}_1} \sgn(A) q^{\|A\|} +
\sum_{A\in\mathcal{G}_2} \sgn(A) q^{\|A\|} \\
&=
\sum_{\substack{A\in\mathcal{G}_1\\  \#A_2 \equiv 0\pmod{6}}} q^{\|A\|} -
\sum_{\substack{A\in\mathcal{G}_1\\ \#A_2 \equiv 3\pmod{6}}}q^{\|A\|} \\
& \qquad +
\sum_{\substack{A\in\mathcal{G}_2\\  \#A_2 \equiv 1\pmod{6}}} q^{\|A\|} -
\sum_{\substack{A\in\mathcal{G}_2\\ \#A_2 \equiv 4\pmod{6}}}q^{\|A\|}.
\end{split}
\end{equation}
As to the left-hand-side of the relation~(\ref{weight-preserving}), we have
\[
\begin{split}
\sum_{A\in\mathcal{F}}\sgn(A) q^{\|A\|}
&=
\sum_{k=-\lfloor\frac{n}{3}\rfloor}^{\lfloor\frac{n+m}{3}\rfloor}
 \sum_{\substack{A\subseteq S\\ \#A=n+3k}} \sgn(A) q^{\|A\|}  \\
&=
 \sum_{k=-\lfloor\frac{n}{3}\rfloor}^{\lfloor\frac{n+m}{3}\rfloor}
 \sum_{\substack{A\subseteq \{1,\ldots,2n+m\} \\ \#A=n+3k}} (-1)^k q^{\|A\|-n(n+3k)} \\
 &=
 \sum_{k=-\lfloor\frac{n}{3}\rfloor}^{\lfloor\frac{n+m}{3}\rfloor}
 (-1)^k q^{-n(n+3k)}  \sum_{\substack{A\subseteq \{1,\ldots,2n+m\} \\ \#A=n+3k}} q^{\|A\|} ,
 \end{split}
 \]
 which with the help of the identity (\ref{weight-sum}) gives
 \begin{equation}\label{lhs-3}
 \begin{split}
 \sum_{A\in\mathcal{F}}\sgn(A) q^{\|A\|}
 &=
 \sum_{k=-\lfloor\frac{n}{3}\rfloor}^{\lfloor\frac{n+m}{3}\rfloor}
 (-1)^k q^{-n(n+3k)} {2n+m\brack n+3k} q^{ {n+3k+1\choose 2}} \\
 &=
 \sum_{k=-\lfloor\frac{n}{2}\rfloor}^{\lfloor\frac{n+m}{2}\rfloor}
 (-1)^k q^{\frac{9 k^2 +3k- n^2+n}{2} }{2n+m\brack n+3k}.
  \end{split}
 \end{equation}
 To evaluate the sum $\sum_{\substack{A\in\mathcal{G}_1\\  \#A_2 \equiv 0\pmod{6}}} q^{\|A\|}$
 in the relation~(\ref{weight-preserving}), we shall  consider partitions into exactly $6l$ distinct parts belonging to the set $\{n,\ldots,n+m\}$ and moreover,
 because of the presence of $0$ among the elements of $\mathcal{S}$, we shall also take into
account partition into exactly $6l-1$ distinct (nonzero) parts belonging to the set $\{n,\ldots,n+m\}$.
 Notice
 that the least sum and the largest sum into exactly $6l$ distinct parts belonging to the set
$\{n, \ldots, n+m \}$ are respectively
\[
n+ (n+1)+ \ldots+  (n+6l-1)=3l(2n+6l-1)
\]
and
\[
(n+m)+(n+m-1)+\ldots+ (n+m-6l+1) = 3l(2n+2m-6l+1)
\]
and similarly the least sum and the largest sum into exactly $6l-1$ distinct parts belonging to the set
$\{n, \ldots, n+m \}$ are  $(6l-1)(n+3l-1)$ and
$(6l-1)(n+m-3l+1)$ respectively.
 Then
 it is easily checked that
\[
 \sum_{\substack{A\in\mathcal{G}_1\\ \#A_2 \equiv 0\pmod{6}}} q^{\|A\|}
=
\prod_{i=1}^{n-1} (q^{-i}+ q^{i}+q^0)
\]
 \begin{equation}\label{rhs-3}
 \begin{split}
& \quad \times \left( \sum_{l=0}^{\lfloor\frac{m+1}{6}\rfloor}
\sum_{j=6l(n+\frac{6l-1}{2})}^{6l(n+m-\frac{6l-1}{2})}
a_{j,6l} q^{j}  +
\sum_{l=0}^{\lfloor\frac{m+2}{6}\rfloor}
\sum_{j=(6l-1)(n+\frac{6l-2}{2})}^{(6l-1)(n+m-\frac{6l-2}{2})} a_{j,6l-1} q^j
\right)  \\
& \qquad\qquad \qquad\qquad\qquad=
q^{\frac{n^2-n}{2}}\frac{(q^3;q^3)_{n-1}}{(q;q)_{n-1}}
A_{6,0}(q).
\end{split}
\end{equation}
Similarly, we obtain
 \begin{align}\label{rhs-3-1}
 \sum_{\substack{A\in\mathcal{G}_1\\  \#A_2 \equiv 3\pmod{6}}} q^{\|A\|}
&=
q^{\frac{n^2-n}{2}} \frac{(q^3;q^3)_{n-1}}{(q;q)_{n-1}}
A_{6,3}(q) ,
\end{align}
 \begin{align}\label{rhs-3-2}
 \sum_{\substack{A\in\mathcal{G}_2\\ \#A_2 \equiv 1\pmod{6}}} q^{\|A\|}
&=
q^{\frac{n^2-n}{2}} \frac{(q^3;q^3)_{n-1}}{(q;q)_{n-1}}
A_{6,1}(q) ,
\end{align}
and
\begin{align}\label{rhs-3-3}
 \sum_{\substack{A\in\mathcal{G}_2\\ \#A_2 \equiv 4\pmod{6}}} q^{\|A\|}
&=
q^{\frac{n^2-n}{2}} \frac{(q^3;q^3)_{n-1}}{(q;q)_{n-1}}
A_{6,4}(q).
\end{align}
Then the desired formula follows by combining the relation (\ref{weight-preserving}) with the formulas (\ref{lhs-3}), (\ref{rhs-3}), (\ref{rhs-3-1}), (\ref{rhs-3-2}), and  (\ref{rhs-3-3}). \\

%%%%%%%%%%%%%%%%%%%%%%%%%%%
\noindent{\bf Acknowledgment.} The author is grateful to the referee
 for valuable comments and suggestions which improved the quality of the paper.
%%%%%%%%%%%%%%%%%%%%%%%%%%%%%%%%%%%%%%%%%%%%%%%%%%%%%%%%%%%%%%%%%%%%%%%%%%%%%

%
\end{document}